# Kendall's tau in high-dimensional genomic parsimony

**Pranab K. Sen**[1]

*University of North Carolina, Chapel Hill*

**Abstract:** High-dimensional data models, often with low sample size, abound in many interdisciplinary studies, genomics and large biological systems being most noteworthy. The conventional assumption of multinormality or linearity of regression may not be plausible for such models which are likely to be statistically complex due to a large number of parameters as well as various underlying restraints. As such, parametric approaches may not be very effective. Anything beyond parametrics, albeit, having increased scope and robustness perspectives, may generally be baffled by the low sample size and hence unable to give reasonable margins of errors. Kendall's tau statistic is exploited in this context with emphasis on dimensional rather than sample size asymptotics. The Chen–Stein theorem has been thoroughly appraised in this study. Applications of these findings in some microarray data models are illustrated.

## Contents

## 1. Introduction

The past three decades have witnessed a phenomenal growth of research literature on statistical methods for large dimensional data models. Such models abound in various interdisciplinary fields, especially in the evolving field of genomics and bioinformatics. *Knowledge discovery and data mining* (KDDM) or statistical learning tools are usually advocated for such high dimensional data models, often on primarily computational or heuristic justifications. The curse of dimensionality is so overwhelming that classical likelihood (principle) based statistical inference tools, baffled with an excessive number of parameters, may not be robust or efficient. Conventional assumptions of multinormality of errors and linearity of regression models

*Supported in part by the C. C. Boshamer Research Foundation at the University of North Carolina, Chapel Hill.

[1]Departments of Biostatistics and Statistics and Operations Research, University of North Carolina, Chapel Hill, NC 27599-7420, USA, e-mail: pksen@bios.unc.edu

*AMS 2000 subject classifications:* Primary 62G10, 62G99; secondary 62P99.

*Keywords and phrases:* bioinformatics, Chen–Stein theorem, dimensional asymptotics, FDR, multiple hypotheses testing, nonparametrics, permutational invariance, $U$-statistics.

may not be generally tenable in such contexts. Moreover, having a large number of coordinate variables, the assumption of their stochastic independence may not be realistic in a majority of cases. On top of that, at least a part of the response variables may be discrete or even purely qualitative in nature; often, the categorical responses may not reveal any (partial) ordering. In that sense, discrete multivariate analysis may appear to be more appropriate than conventional multinormal model based analysis. Even for multinormal models, the high-dimensionality may demand a far larger sample size in order to implement a full likelihood based asymptotic analysis. That is, we need the conventional $n \gg K$ environment for drawing appropriate statistical conclusions with reasonable precision.

Typically, in such high-dimensional models, one encounters a $K \gg n$ environment, where $K$ is the dimension of the data and $n$ is the sample size. In such *high-diensional low sample size*, HDLSS, models, effective dimension reduction may be a challenging statistical task, usually beyond the scope of KDDM. For example, in neuronal spike train models, there are literally tens of thousands of neurons (nerve cells), and in the presence of external stimuli, the spike trains for any observable subset of neurons exhibit a high-degree of nonstationarity. Further, recording of such spike trains in a large number of nerve cells may be invasive to the brain functioning due to the destructive nature of recording ([16], Ch. 3). Essentially, we have a very high dimensional counting process. Doubly stochastic Poisson processes have been considered in the literature, albeit without much claim of optimal resolutions. In magnetic resonance imaging, MRI, there could be tens of thousands of microscopic units producing an enormously high dimensional spatial data model. More complexities may arise in case of (functional) fMRI models. For such HDLSS models, parametric asymptotics may not have adequate scope or good statistical interpretation.

The transition from conventional normal theory to nonparametric linear models has been well fortified along with the development of nonparametric or robust statistical methods based on $R$- statistics (ranks), $M$-statistics (maximization) and linear combinations of order statistics or $L$-statistics; see, for example, [8] where other pertinent references have been extensively cited. In a more general setup, nonparametric regression functionals have been formulated wherein the linearity of regression or a specific nonlinear form are not assumed to hold.

In the context of testing monotonicity of nonparametric regression, without assuming a linear or any specific nonlinear form, Ghosal et al. [5] considered suitable $U$-processes based on a locally smoothed Kendall's tau statistic. They provided general asymptotics for such locally smoothed Kendall's tau processes when both the independent and dependent variates are stochastic, and illustrated their effective use in the postulated hypothesis testing problem. Such local versions of Kendall's tau statistics have simple statistical interpretation, albeit, in view of possibly slower rate of convergence, the impact of large sample size is apparent in their analysis. In the contemplated bioinformatics area, as we shall see, the HDLSS scenario calls for alternative approaches, and some of these will be explored in this study.

In a simple regression setup, the Theil-Sen (point as well as interval) estimates of the regression slope based on the Kendall tau statistic [15], have simple forms, and are computationally tractable and statistically robust. Another advantage of the Kendall tau statistic is its adaptability for count data as well as latent-effect models. Further, a test for the null hypothesis of no regression based on the Kendall tau statistic (being distribution-free under the null hypothesis of invariance) remains valid and efficient for such complex models. Our contemplated models, unlike [5], entail a high dimensional data with relatively (and often inadequately) smaller sample size, i.e., the HDLSS ($K \gg n$) environment. As we shall see in the next section, there may not be a genuine temporal pattern. In addition, there may be other complications arising from lack of spatial-compactness, spatial homogeneity and other spatial dependence patterns.

For better motivation, in Section 2, an illustration is made with a microarray data model where HDLSS models typically arise. Section 3 deals with the appropriateness of statistical modeling and analysis based on a pseudo-marginal approach incorporating coordinatewise construction of the Kendall tau statistic, in such $K \gg n$ environments. Section 4 is devoted to the dimensional asymptotics for the Kendall tau process in such HDLSS models where there are two basic problems : (i) group divergence, and (ii) classification of genes into disease and nondisease types. For the first problem, a pseudo-marginal approach based on the Hamming distance has been explored in [18] while in the latter context, multiple hypotheses testing (MHT) problems in HDLSS setups arise in a different perspective and call for some alternative novel tools for valid and efficient statistical appraisals. Motivated by these perspectives in such HDLSS models, some applications of the Chen–Stein [3] theorem in such $K \gg n$ environments are presented in the last section. These generalizations cover both the MHT and the gene-environment interaction testing problems.

## 2. An illustrative data model

We consider a genomic model arising in microarray data analysis as an illustration. The microarray technology allows simultaneous studies of thousands of genes, $K$, possibly differentially expressed under diverse biological/experimental setups, with only a few, $n$, arrays. We may refer to Lobenhofer et al. [11] where for a set of 1900 genes, arranged in rows, the gene expressions were recorded at 6 time points, with 8 observations at each time point. Thus $1900 = K \gg n = 48$. The gene-expression levels are measured by their color intensity (or luminosity) as a quantitative (nonnegative) variable, either on the (0, 1) or 0–100 per cent scale, or (based on the log-scale) on the real line $\Re$. A gene associated (causally or statistically) with a target disease is known as a *disease gene*, DG, while the others as *nondisease genes*, NDG. Gene expression levels under different environments cast light on plausible *gene-environment interactions* (or associations) so that if the arrays are properly designed, *mapping disease genes* may be facilitated with such microarray studies. One of the main issues is identifying differentially expressed genes among thousands of genes, tested simultaneously, across experimental conditions. Typically, for a target disease, there are only a few DG while the NDG comprise the vast majority. A NDG is expected to have a low gene expression level while a DG is expected to have generally higher expression levels. Thus, a natural *stochastic ordering* of gene expression levels of the DG with varying disease severity is plausible while the NDG expression levels are expected to be stochastically unaffected by such disease level differentials.

Microarray data go thorough a lot of standardization and normalization so that conventional simple models, such as the classical MANOVA models, may rarely be totally adaptable. If the arrays are indexed by an explanatory or design variate ($t$) that possesses an ordering (not necessarily linear), then the stochastic ordering could be exploited through suitable nonparametric techniques. The main difficulty in modeling and statistically analyzing microarray data stems from the high dimensionality of the genes compared to the number of arrays. While the different arrays may sometimes be taken to be at least statistically independent, the genes may not. Moreover, not much is known about the spatial topology of the genes or their genetic distances. There is another factor that merits our attention. The gene expression levels for the different genes in an array are neither expected to be stochastically independent nor (marginally) identically distributed. Sans such an i.i.d. clause, standard parametrics typically adaptable for fMRI models (albeit mostly done in a Bayesian coating) may encounter roadblocks for fruitful adaptation in microarray data models. Thus, structurally, such data models are different from those usually encountered in nonparametric functional regression models. For this reason, a pseudo-marginal approach is highlighted here. This approach exploits the marginal nonparametrics fully and renders some useful modeling and analysis convenience.

## 3. Some HDLSS formulations

Motivated by microarray data models introduced in Section 2, we consider here a set of $n$ arrays (sample observations) where there is a design variate $t_i$ associated with the $i$th array, for $i = 1, \ldots, n$. Without loss of generality, we assume the $t_i$ are ordered, i.e.,

$$(3.1) \qquad t_1 \leq t_2 \leq \cdots \leq t_n,$$

with at least one strict inequality. We do not, however, impose any linear or specific parametric ordering of these design variates. The multisample (ordered alternative) model is a particular case where $n$ can be partitioned into $I$ subsets of sizes $n_1, \ldots, n_I$ such that within each subgroup, the $t_i$ are the same while they are ordered over the $I$ different subsets. For the $i$th array, corresponding to the $K$ genes (positions), we have a gene expression level denoted by $X_{ik}$, $k = 1, \ldots, K$, so that we have $K$-vectors $\mathbf{X}_i = (X_{i1}, \ldots, X_{iK})'$, for $i = 1, \ldots, n$. The joint distribution function of $\mathbf{X}_i$ is denoted by $F_i(\mathbf{x})$, $\mathbf{x} \in \Re^K$. Further, for the $k$th gene in the $i$th array, i.e., $X_{ik}$, the marginal distribution is denoted by $F_{ik}(x)$, $x \in \Re$, for $k = 1, \ldots, K$; $i = 1, \ldots, n$. For a given $i$, the $F_{ik}$, $k = 1, \ldots, K$ may not be generally the same, and moreover, the $X_{ik}$, $k = 1, \ldots, K$ may not be all stochastically independent.

If a gene $k$ is NDG and the $t_i$ reflect the variability of the disease level, then the $F_{ik}, i = 1, \ldots, n$ should be the same. On the other hand, for a DG $k$, for $i < i'$, $X_{ik}$ should be stochastically smaller than $X_{i'k}$ in the sense that the $F_{ik}, i = 1, \ldots, n$ should have the ordering

$$(3.2) \qquad F_{1k}(x) \geq F_{2k}(x) \geq \cdots \geq F_{nk}(x), \forall\ x \in \Re.$$

Therefore, we could force a characterization of DG and NDG based on the following stochastic ordering: For a NDG $k$, the $F_{ik}$, $i = 1, \ldots, n$ are all the same, this being denoted by the null hypothesis $H_{0k}$, while for a DG $k$, the stochastic ordering in (3.2) holds which we denote by $H_{1k}$, for $k = 1, \ldots, K$. In this marginal formulation, we have a set of $K$ hypotheses corresponding to the $K$ genes, and whatever appropriate test statistic (say $T_{nk}$) we use for testing $H_{0k}$ vs. $H_{1k}$, these statistics may not be, generally, stochastically independent. The basic problem is therefore to test simultaneously for

$$(3.3) \qquad H_0 = \bigcap_{k=1}^{K} H_{0k} \text{ vs } H_1 = \bigcup_{k=1}^{K} H_{1k},$$

without ignoring possible dependence of the test statistics for the component hypotheses testing $H_{0k}$ vs $H_{1k}$, for $k = 1, \ldots, K$. This makes it appealing to follow the general guidelines of the Roy [13] *union-intersection principle* (UIP), albeit in a marginalization (i.e., adapting a finite union and finite intersection scheme), and thus permitting a more general framework so as to allow simultaneous testing and classification into DG / NDG groups. In a very parametric setup, some order restricted inference problems have been considered by [12]. However, in our setup, such normality based parametric models may not be very appropriate.

Our approach is based on the classical Kendall tau statistics for each of the $K$ genes and the incorporation of these (possibly dependent) marginal statistics in a composite scheme for classification. For the $k$th gene, based on the $n$ observations $X_{ik}$, $i = 1, \ldots, n$, and the tagging variables $t_1, \ldots, t_n$, we define the Kendall tau statistic as

$$(3.4) \qquad T_{nk} = \binom{n}{2}^{-1} \sum_{1 \leq i < i' \leq n} \operatorname{sign}(X_{i'k} - X_{ik})\operatorname{sign}(t_{i'} - t_i),$$

for $k = 1, \ldots, K$. Conventionally, we take $\operatorname{sign}(0) = 0$. Note that $T_{nk}$ is a (generalised) $U$-statistic of degree 2 [7]. Further, note that by (3.1), we may set $\mathcal{S} = \{(i, i') : t_i < t_{i'}; 1 \leq i < i' \leq n\}$ and let $N$ be the cardinality of the set $\mathcal{S}$. Then by (3.1), $n - 1 \leq N \leq \binom{n}{2}$. Moreover, we may rewrite $T_{nk}$ as

$$(3.5) \qquad T_{nk} = \binom{n}{2}^{-1} \sum_{\mathcal{S}} \operatorname{sign}(X_{i'k} - X_{ik}), \ k = 1, \ldots, K,$$

where $\mathcal{S}$ depends on the ordering of the $t_j$ and therefore remains the same for every $k = 1, \ldots, K$. Note further that whenever $N < \binom{n}{2}$, the range of variation of $T_{nk}$ is $\left(-N/\binom{n}{2}, N/\binom{n}{2}\right)$ which is contained in the interval $(-1, 1)$. That is why we shall find it convenient to take the modified or rescaled Kendall tau as

$$(3.6) \qquad T_{nk}^o = N^{-1} \sum_{\mathcal{S}} \operatorname{sign}(X_{i'k} - X_{ik}),$$

whose range is exactly $(-1, 1)$, albeit the distribution being still discrete.

Note that for any $k = 1, \ldots, K$, under $H_{0k}$, for every $i \neq i'$, the difference $X_{i'k} - X_{ik}$ is symmetrically distributed around 0, and hence, $E_{0k}\{\operatorname{sign}(X_{i'k} - X_{ik})\} = 0$ so that

$$(3.7) \qquad E_{0k}\{T_{nk}\} = E_{0k}\{T_{nk}^o\} = 0, \ \forall \ k = 1, \ldots, K.$$

Further, the marginal distribution of $T_{nk}$ under $H_{0k}$ is generated by the $n!$ equally likely permutations of the $X_{ik}$ among themselves. Therefore when all the $F_{ik}$ are continuous, ties among the observations being negligible with probability 1, $T_{nk}$ (or $T_{nk}^o$) is distribution-free under $H_{0k}$. This distribution may depend on the set $\mathcal{S}$ but that being the same for all $k$, we conclude that under $H_0$ in (3.3), marginally each $T_{nk}^o$ is distribution-free and these $K$ statistics all have the same marginal distribution. If all the $t_i$ were stochastic and continuous then $N = \binom{n}{2}$ and we will have

$$(3.8) \qquad \operatorname{Var}_0\{T_{nk}\} = 2(2n + 5)/\{9n(n - 1)\}.$$

On the other hand, in general for $N \leq \binom{n}{2}$, the $t_i$ are not distinct and may be even nonstochastic, and hence, the variance is equal to

$$(3.9) \qquad \operatorname{Var}_0\{T_{nk}^o\} = N^{-2}\{(2/3)(N_1 - N_2) + N\} = \nu_n^2, \text{ say,}$$

where $N_1$ is the cardinality of the set $\{(i, i'), (i, i'') : t_i < t_{i'}, t_i < t_{i''}, t_{i'} \neq t_{i''}\}$ and $N_2$ is the cardinality of the set $\{(i, i'), (i'', i) : t_i > t_{i''} \neq t_{i'}\}$. For small values of $n$ and given (3.1), one can enumerate $\mathcal{S}$ and obtain the exact distribution of $T^o_{nk}$ under $H_{0k}$. If $n$ is large, the standardized form of the statistic, i.e., $T^o_{nk}/\nu_n$ has closely a standard normal distribution. In our setup, perhaps the exact permutation distribution plays a greater role and this will be illustrated later on.

The behavior of $T^o_{nk}$ under alternatives would naturally depend on the stochastic ordering in (3.2) and these statistics will not be exact distribution-free nor possibly have identical marginal laws. Nevertheless, under (3.2), for every $i < i'$, $X_{i'k} - X_{ik}$ has a distribution tilted to the right, so that

$$(3.10) \qquad E\{T^o_{nk} \mid H_{1k}\} \geq 0, \ \forall \ k = 1, \ldots, K.$$

This motivates us to use tests based on the marginal statistics $T^o_{nk}$ using the right hand side critical region, or equivalently the right-hand sided $p$-values. Recall that the distribution of each $T^o_{nk}$, at least for $n$ not too large, is discrete, but that is not going to be of any particular concern. A greater concern is to incorporate possible stochastic dependence among the $K$ statistics $T^o_{nk}$, $k = 1, \ldots, K$ (even under the null hypothesis) and their possible heterogeneity when some of the $H_{1k}$ are true. A basic problem is to formulate suitable multiple hypothesis testing procedures to assess which hypotheses are to be rejected subject to a suitably defined Type I error rate. This is elaborated in the next section.

## 4. Dimensional asymptotics and the union intersection test

Although independence across microarrays may be assumed, their i.d. structure may be vitiated if the arrays relate to different biological or experimental setups. Moreover, for different genes, the gene expression (marginal) distributions are likely to be different when there is gene-environment interaction. Taking into account such plausible inter-gene stochastic dependence and heterogeneity, we need to prescribe statistical modeling and analysis tools. This will be accomplished through dimensional asymptotics where $K$ is made to increase indefinitely while $n$, being small compared to $K$, may or may not be adequately large.

In view of (3.3), it is tempting to appeal to the union-intersection principle [13], or UIP, to construct suitable test statistics which will cover the genome-wise picture in a reasonable way. Towards this, we may note that as under $H_0$ (i.e., $H_{0k}, \forall k$), marginally each $T^o_{nk}$ has the same distribution (which does not depend on the underlying $F_{ik}$). Thus, corresponding to any $c : -1 \leq c \leq 1$, the tail probability $P_0\{T^o_{nk} > c\}$ is the same for all $k$ and this can be evaluated by using the exact permutation distribution generated by the $n!$ permutations of the $X_{ik}, 1 \leq i \leq n$. The UIP then leads to the following union-intersection test, UIT, statistic:

$$(4.1) \qquad T^{*o}_n = \max\{T^o_{nk} : 1 \leq k \leq K\},$$

where the test function is given by $\phi(T^{*0}_n) = 1, \gamma$, or 0, accordingly as $T^{*o}_n$ is $>, =$ or $< c$ and $\gamma : (0 \leq \gamma \leq 1)$ is so chosen that $E_0\{\phi(T^{*o}_n)\} = \alpha$, the preassigned *level of significance*. Note that for $n$ not adequately large, the null distribution of $T^o_{nk}$ is essentially discrete and hence this usual randomization test function is aimed to take care of this problem.

The crux of the problem is therefore to determine such a critical level $c_\alpha$. The joint distribution of the $T^o_{nk}, 1 \leq k \leq K$, even under the null hypothesis $H_0$, depends on the underlying $K$-dimensional distribution $F_i$, and hence, in general will not be distribution-free. Thus, the usual technique of finding out the critical level of $T_n^{*o}$ from this joint distribution may be intractable.

One possibility is to incorporate the fact that under $H_0$, the $K$-vectors $\mathbf{X}_i, i = 1, \ldots, n$, are i.i.d. and hence their joint distribution remains invariant under any permutation of these vectors among themselves. Thereby we can evaluate such critical values by an to appeal to the permutation distribution generated by the $n!$ equally likely permutations of the $K$-vectors $\{\mathbf{X}_i\}$ among themselves. This permutation law generates the (unconditional null) marginal laws of the $T_{nk}^o$, and provides some conditional versions of their joint distributions of various orders. Since this permutation law is a conditional law (given the collection of all these $K$-vectors), the critical values obtained in this manner are themselves stochastic, thus introducing another layer of variation. Nevertheless, it provides a conditionally distribution-free test. One discouraging feature of this permutation approach is that the permutation invariance does not hold under the alternative hypothesis, and hence critical levels computed from the permutation law involving an observed set of $\{\mathbf{X}_i\}$ may be sensitive to the data conformity to the null situation.

If we assume that all the $T_{nk}^o$ are stochastically independent, then we have for any $c$, $-1 \leq c \leq 1$, under $H_0$,

$$(4.2) \qquad P_0\{T_n^{*o} \leq c\} = [P_0\{T_{n1}^o \leq c\}]^K,$$

so that the distribution-free nature of the $T_{nk}$ under the null hypothesis provides the access to the computation of the test function and the critical level. If $n$ is at least moderately large, in view of the asymptotic normality of $T_{nk}^o/\nu_n$, the randomization test function may be replaced by a conventional normal theory test function, where for the individual tests, a significance level $\alpha^*$ is so chosen that

$$(4.3) \qquad \alpha = 1 - (1 - \alpha^*)^K.$$

Generally, if we let $\alpha^* = (\alpha/K)$, then the size of the UIT is $\leq \alpha$ no matter whether the $T_{nk}^o$ are stochastically independent or not. There is, therefore, a certain amount of conservativeness in this specification.

In passing, we may remark that by the classical asymptotics on Hoeffding's $U$-statistics, any pair $(T_{nk}^o, T_{nq}^o)$, with $k \neq q$, is a bivariate $U$-statistic, for $\alpha^*$ sufficiently small, so using the bivariate extreme statistics results (viz., [19]), we can claim that the events $\{T_{nk}^o > c_{\alpha^*}\}$ and $\{T_{nq}^o > c_{\alpha^*}\}$ will be asymptotically (as $K \to \infty$) independent so that $P_0\{T_{nk}^o > c_{\alpha^*}, T_{nq}^o > c_{\alpha^*}\}$ can be well approximated by $[P_0\{T_{nk}^o > c_{\alpha^*}\}]^2$. In a similar manner, the third order probability terms can be handled, and the Bonferroni bound retaining the second and third order probabilities provide a good approximation : $\alpha = K\alpha^* - \binom{K}{2}\alpha^{*2} + \binom{K}{3}\alpha^{*3} + o(\alpha^{*3})$. As a result, $\alpha^* = (\alpha/K)$ provides a good approximation to the level of significance. Therefore, for the UIT, when $K$ is large, even when the genes are not stochastically independent, letting $\alpha^* = (\alpha/K)$ we may consider the following multiple hypothesis testing scheme:

*For a chosen $\alpha^* = K^{-1}\alpha$, obtain the marginal distributional critical level $c_{\alpha^*}$, and reject those $H_{0k}$; $k \in \{1, \ldots, K\}$ for which the corresponding $T_{nk}^o$ exceeds $c_{\alpha^*}$.*

A randomization test function can be prescribed when $n$ is not adequately large. Thus, the UIT provides a bound on the *family wise error rate*, FWER. If we take $\alpha* \sim \alpha/K$ and $K$ is large, we need to make sure that $n$ is so large that $\nu_n^{-1} c_{\alpha^*} < 1$; this will imply that if we are to use the permutation null distribution of any $T^o_{nk}$, being attracted by the permutational central limit theorem, it has a nonzero mass point beyond $c_{\alpha^*}/\nu_n$. If $\nu_n^2 = O(n^{-1})$, as is typically the case, then $c_{\alpha^*} = O(n^{-1/2}\sqrt{-2\log \alpha^*})$ so that $\log K = O(n)$ and this does not appear to be a serious concern in real life applications. For example, if we have three groups of arrays, say within each group there are 5 arrays, the total number of partitioning 15 units into 3 subsets of 5 each is equal to $(15)!/(5!)^3$ and this is so large (756,756) that even if $K$ is as large as 30,000, it would not be a problem. However, for large $K$, the UIT, like the classical likelihood ratio test, will have little power, and hence alternative test procedures need to be explored. This illustrates the important role of the design of the study and the number of arrays required in trying to include a very large $K$.

Roy's UIT can be adapted by exploring the information contained in the ordered $p$-values. If the $T^o_{nk}$ are all stochastically independent (and as they are identically distributed under the null hypothesis $H_0$) then one can adapt Simes' [20] theorem (which is a restatement of the classical Ballot theorem (viz., [9]) introduced some twenty years earlier). If $P_1, \ldots, P_K$ are the $p$-values for the $K$ marginal tests and $P_{K:1} \leq \cdots \leq P_{K:K}$ are the corresponding order statistics, then assuming that under $H_0$ the $P_k$ have a uniform $(0, 1)$ distribution (i.e., tacitly assuming that the $T^o_{nk}/\nu_n$ have a continuous distribution under $H_0$), Simes' theorem asserts that for every $\alpha : 0 < \alpha < 1$,

$$(4.4) \qquad P\{P_{K:k} > k\alpha/K, \ \forall \ k = 1, \ldots, K \ |H_0\} = 1 - \alpha.$$

Suppose now we define the *anti-ranks* $S_1, \ldots, S_K$ by letting

$$(4.5) \qquad P_{K:k} = P_{S_k}, \ k = 1, \ldots, K,$$

where again ties among the ranks are neglected under the assumption of continuity of the distribution of the $P_k$. Whereas Simes' theorem provides a test of the overall hypothesis, Hochberg [6] derived a step-up procedure for multiple hypotheses testing based on the following : For every $\alpha \in (0, 1)$,

$$(4.6) \qquad P\{P_{K:k} \geq \alpha/(K - k + 1), \ \forall \ k = 1, \ldots, K \ |H_0\} = 1 - \alpha.$$

Benjamini and Hochberg [2] considered a step-up procedure based on the Simes theorem. Their multiple hypothesis testing procedure is the following:

*Reject those null hypotheses* $\{H_{0S_k}\}$ *for which* $P_{S_k} \leq k\alpha/K, \ k = 1, \ldots, K$, *and accept those null hypotheses in the complementary set.*

For some related developments in a parametric setup, we refer to [2], [4], [10], [14] and [21], among others.

These developments paved the way for other measures of error rates which are more adaptable in the $K \gg n$ environment. Some of these will be discussed later on. There are two basic concerns that can be voiced in this respect. The whole setup is based on the assumed uniform distribution of the $P_k$ under the null hypothesis. However, if we look into the statistics $T^o_{nk}$ in our setup, we may note that though they have a specified distribution, the latter is a discrete one defined over the interval $(-1, 1)$. Noting that there are a set of discrete mass points, ties among the $T^o_{nk}/\nu_n$ (and hence $P_k$) can not be neglected with probability one, and moreover, the $P_k$ will have a set of probability mass points on $[0, 1]$ with non-zero masses. Thus, technically the above probability results are not strictly usable (unless $n$ is indefinitely large, contradicting the $K \gg n$ environment). Secondly, as was stressed earlier, the $T^o_{nk}$ across the set of genes are generally not stochastically independent.

Controlling the FWER when $K$ is very large may generally entail undue conservativeness of multiple hypotheses testing schemes. On the other hand, using a level of significance for each marginal hypothesis testing problem may lead to a large FWER.

In the context of microarrays suppose that there are $K_1$ disease genes (DG) and $K_0 = K - K_1$ NDG; thus, we have a set of $K_0$ null hypotheses which are true and a complementary set of $K_1$ hypotheses which are not true. Suppose that based on our multiple hypotheses testing procedure, we accept $m_0$ out of $K_0$ true null hypothesis so that the remaining $K_0 - m_0 = m_1$ true null hypotheses are rejected. Similarly, among the $K_1$ not true null hypotheses, $l_0$ are accepted as true and $l_1$ accepted in favor of the alternative. Thus, a totality of $R = m_1 + l_1$ hypotheses are rejected while $K - R$ are accepted. Mind that though we observe $R$, through our chosen multiple hypotheses testing procedure, individually $m_1, l_1$ are not observable; all these $(R, l_1, m_1)$ are stochastic in nature. A natural modification of the FWER, to suit such $K \gg n$ environments, is the *per-comparison error rate* (PCER) defined as

$$(4.7) \qquad \text{PCER} = E(m_1)/K,$$

which is the expected proportion of Type I errors among the $K$ hypotheses. A related measure is the *per-family error rate* (PFER), defined as

$$(4.8) \qquad \text{PFER} = E(m_1),$$

which is the expected total number of Type I errors among the $K$ hypotheses. Obviously, PFER $= K.$PCER, and is generally large when $K$ is large (unless the PCER is very small). Moreover,

$$(4.9) \qquad \text{PFER} = E(m_1) = \sum_{r \geq 1} rP\{m_1 = r\} \geq P\{m_1 > 0\},$$

so that PFER $\geq$ FWER.

If our observed $R = 0$ then no true null hypothesis is rejected and hence there is no false discovery. For $R \geq 1$, the proportion of false discovery is given by $Q = m_1/R$; conventionally, it is taken $Q = 0$ when $R = 0$, so that $Q$ is properly defined for every nonnegative $R$ and $m_1$. However, $Q$ is not observable. Hence, the *false discovery rate* (FDR) is defined as

$$(4.10) \qquad \text{FDR} = E\{Q\} = \sum_{r \geq 1} P\{R = r\} E\{m_1/R | R = r\}.$$

Since, conventionally, we have forced $Q = 0$ for $R = 0$, this definition of FDR may produce a negative bias. An alternative definition, known as the $^pFDR$, is defined as

$$(4.11) \qquad ^p\text{FDR} = E\{Q|R > 0\} = \text{FDR}/P\{R \geq 1\}.$$

Naturally, $^p$FDR $\geq$ FDR.

In the formulation of FDR and $^p$FDR it is not necessary to assume that all of the test statistics have continuous distributions under the null hypothesis. If these distributions are all continuous then of course the $p$-values have a uniform $(0, 1)$ distribution under the null hypothesis, and hence, the multiple hypotheses testing schemes discussed earlier can be conveniently adapted. In our setup, each test statistic has marginally the same null distribution, albeit that is discrete. So, it might be necessary, especially when $n$ is not large, to make use of this otherwise completely specified, discrete distribution without assuming a uniform distribution for the associated $p$-values under the null hypothesis.

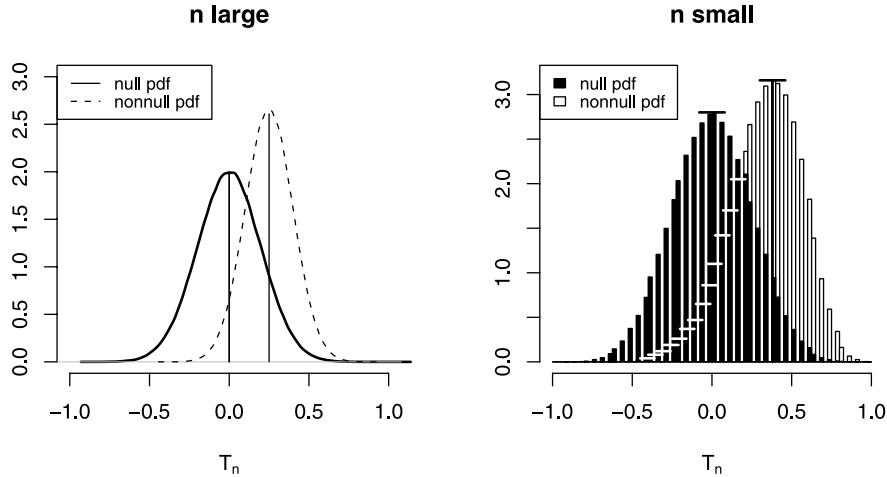

Fig 1. *Comparison of the null distribution with the alternative distribution.*

We may simulate the permutation distribution of any marginal test statistics and thereby take into account possible dependence among the gene expressions without assuming any specific pattern. Of course, marginally, each test statistic has the same null distribution. So, if we consider the set $\{T_{nk}^o : k = 1, \ldots, K\}$ and define the empirical distribution

$$(4.12) \qquad G_K(t) = K^{-1} \sum_{k=1}^{K} I(T_{nk}^o \leq t), \ t \in (-1, 1),$$

then $E_0\{G_K(t)\} = G(t), \forall t \in (-1, 1)$ where $G(t)$ is the common marginal distribution of the $T_{nk}^o$ under the null hypothesis. The summands in $G_K(t)$ are all bounded variables, nondecreasing in $t \in (-1, 1)$ and $G(t)$ is also nondecreasing and assumes values on $(0, 1)$. Thus, whenever $G_K(t)$ stochastically converges pointwise to $G(t)$, it does so uniformly in $t \in (-1, 1)$. Further $G_K(t) - G(t)$ is a bounded r.v., and hence, if it converges in probability, it converges in the $r$th mean for every $r > 0$. Therefore it might suffice to assume that the dependence pattern satisfies the condition:

$$(4.13) \qquad \operatorname{Var}(G_K(t)) \to 0, \text{ as } K \to \infty.$$

Then we conclude that $\|G_K(.) - G(.)\| = \sup\{|G_K(t) - G(t)| : t \in (-1, 1)\}$ stochastically converges to 0. Further, (4.13) holds under quite general dependence patterns.

It is naturally tempting to explore weak convergence (invariance principles) results for $\sqrt{K}(G_K(.) - G(.))$ wherein $K$ is taken indefinitely large but not $n$. Since $G(t), t \in (-1, 1)$ is a discrete distribution function with mass points over $(-1, 1)$, the jump-discontinuities of $G(.)$ may vitiate the usual compactness (or tightness) properties possessed in the continuous case, albeit by strengthening (4.13) to

$$(4.14) \qquad \limsup_K K \operatorname{Var}(G_K(t)) < \infty, \forall \ t \in (-1, 1),$$

pointwise, the asymptotic normality (as $K \to \infty$) follows under quite general dependency conditions. If we have some linear functional of $G_K(.)$ as a test statistic, this weak convergence would have been quite useful in deriving the asymptotic (in $K$) normality of the test statistic under the null hypothesis; (4.14) would have been sufficient in that context. However, in our case, we have some functional of $G_K(.)$, of extremal order statistic type, namely, the extreme quantiles of a set of dependent r.v.s, and hence we may need somewhat different regularity conditions. This perspective is appraised more elaborately in the next section.

## 5. Dimensional asymptotics and Chen–Stein theorem

In the previous section we have briefly discussed the plausibility of some $K_o$ NDG and $K_1$ DG with $K_o + K_1 = K$, the total number of genes. Neither $K_1$ nor the DG positions are known and hence we have a dual problem of estimating $K_1$ as well as identifying the positions of these $K_1$ DG's. It is conceivable that the NDG having stochastically smaller expression levels (than the DG) and the stochastic dependence among the DG may not be insignificant. We intend to incorporate this stochastic dependence structure among the gene expressions in a suitable model. Unfortunately, sans any positional ordering of the $K$ genes, it might be difficult to assume suitable mixing conditions under which central limit theorems may apply. As for considering alternative limit theorems for dependent sequences, we intend to incorporate the Chen–Stein theorem [3] and its ramifications wherein Poisson approximations for more general dependent sequences are advocated. For our convenience, let us state the Chen–Stein Theorem in a slightly updated version [1].

**Theorem 1.** *(Chen–Stein): Let $\mathcal{I}$ be an index set with elements $i \in \mathcal{I}$ and let $K$ be the cardinality of the set $\mathcal{I}$. For each $i \in \mathcal{I}$ let $Y_i$ be an indicator random variable and let*

(5.1) $$P\{Y_i = 1\} = 1 - P\{Y_i = 0\} = p_{Ki}, \ i \in \mathcal{I}.$$

*Let $W = \sum_{i \in \mathcal{I}} Y_i$ the total number of occurrence of the events $\{Y_i = 1\}$, $i \in \mathcal{I}$, and let $\lambda_K = \sum_{i \in \mathcal{I}} p_{Ki} = E(W)$. For each $i \in \mathcal{I}$, we define a set $\mathcal{J}_i \in \mathcal{I}$ and its complement $\mathcal{J}_i^c$ as the set of dependence of $i$ and its complement, set of independence of $i$. Thus, it is tacitly assumed that $Y_i$ is independent of $\{Y_j, j \in \mathcal{J}_i^c\}$, for every $i \in \mathcal{I}$. Further, let*

$$b_1 = \sum_{i \in \mathcal{I}} \sum_{j \in \mathcal{J}_i} E(Y_i)E(Y_j);$$

(5.2) $$= \sum_{i \in \mathcal{I}} \sum_{j \in \mathcal{J}_i} p_{Ki} p_{Kj},$$

(5.3) $$b_2 = \sum_{i \in \mathcal{I}} \sum_{j(\neq i) \in \mathcal{J}_i} E(Y_i Y_j),$$

*and*

(5.4) $$b_3 = \sum_{i \in I} E|\{E(Y_i - E(Y_i))|\{Y_j, \forall j \in \mathcal{J}_i^c\})|.$$

Finally, let $Z$ be a random variable having Poisson distribution with parameter $E(Z) = \lambda_K$. Then

$$\|\mathcal{L}(W) - \mathcal{L}(Z)\| \leq 2(b_1 + b_2 + b_3)\frac{1 - e^{-\lambda_K}}{\lambda_K}$$
(5.5)
$$\leq 2(b_1 + b_2 + b_3)\min\{1, \lambda_K^{-1}\}.$$

A direct corollary to Theorem 1 is the following:

(5.6) $$|P\{W = 0\} - e^{-\lambda_K}| \leq 2(b_1 + b_2 + b_3)\min\{1, \lambda_K^{-1}\}.$$

An interesting feature of this Theorem is the dual control of $\lambda_K$, the expectation and $b_1, b_2$, and $b_3$, the dependence functions. In line with our intended application we consider a natural extension of this result. With the same notation as in Theorem 1, we replace the $Y_i, i \in \mathcal{I}$, by a sequence of processes $Y_i(t), i \in \mathcal{I}, t \in T$, where $T = (0, a)$, for some $a > 0$, and assume that for each $i$, $Y_i(t)$ is nondecreasing in $t$ and yet a zero-one valued random variable. Further assume that the sets $\mathcal{J}_i$ do not depend on $t \in T$. For every $i \in \mathcal{I}$, $t \in T$, we denote by $p_{Ki}(t) = E(Y_i(t))$, and the corresponding parameters by $\lambda_K(t)$, $b_1(t), b_2(t)$ and $b_3(t)$. Let $\mathbf{W}_K = \{W_K(t), t \in T\}$ be the sum process and corresponding to $Z$, we introduce a Poisson process $\mathbf{Z}_K = \{Z_K(t), t \in T\}$ whose expectation process is $\{\boldsymbol{\lambda}_K = \{\lambda_K(t), t \in T\}$. Then

$$\|\mathcal{L}(\mathbf{W}_K) - \mathcal{L}(\mathbf{Z}_K)\| \leq 2\sup\{(b_1(t) + b_2(t) + b_3(t))\frac{1 - e^{-\lambda_K(t)}}{\lambda_K(t)} : t \in T\}.$$

The proof of this extension is along the lines of Theorem 1 and hence we omit the details.

In our study, unless $n$ is large, we may not have a continuous time parameter ($t \in T$). Thus, we consider an intermediate result that remains applicable for small $n$ as well.

**Theorem 2.** *Consider a set of $M$ discrete time points $-1 \leq \tau_1 < \cdots < \tau_M \leq 1$ with respective probability masses $\eta_{n1}, \ldots, \eta_{nM}$ where $M$ may depend on $n$. Also, let $\nu_{nj} = \sum_{i \leq j} \eta_{ni}$, $j = 1, \ldots, M$. Further, let $Y_i(\tau_j), i = 1, \ldots K$, $j = 1, \ldots, M$ be an array of zero-one valued random variables where $Y_i(\tau_j)$ is nondecreasing in $\tau_j$ and $E(Y_i(\tau_j)) = \nu_{nj}$, $j = 1, \ldots, M$. Define $\mathbf{W}_K = \{W_K(\tau_j), j = 1, \ldots, M\}$ where $W_K(\tau_j) = \sum_{i=1}^{K} Y_i(\tau_j)$ for $j = 1, \ldots, M$. Similarly, let $\mathbf{Z}_K = \{Z_K(\tau_j), j = 1, \ldots, M\}$ be a discrete time parameter Poisson process with the drift function $\boldsymbol{\nu}_K = \{\nu_{nj}, j = 1, \ldots, M\}$. Define the parameters $b_{K1}(\tau_j), b_{K2}(\tau_j), b_{K3}(\tau_j)$, $j = 1, \ldots, M$ as in (5.2), (5.3), and (5.4); assume that as $K \to \infty$,*

(5.7) $$\max\{(b_{K1}(\tau_j) + b_{K2}(\tau_j) + b_{K3}(\tau_j))\frac{1 - e^{-\nu_{nj}}}{\nu_{nj}} : j \leq M\} \to 0.$$

*Then, as $K$ increases indefinitely,*

(5.8) $$\|\mathcal{L}(\mathbf{W}_K) - \mathcal{L}(\mathbf{Z}_K)\| \to 0.$$

Again, being a finite-dimensional version of Theorem 1, this does not need an elaborate proof.

In the present context, under the null hypothesis, all the $T^o_{nk}$ have a common distribution over $(-1, 1)$; this is discrete but symmetric about 0, and is completely known (though could be computationally intensive if $n$ is not too small). Let us denote the distinct mass points for $T_{nk}^o$ by $-1 = a_1 < a_2 < \cdots, a_L = 1$ and let

(5.9) $$\tau_j = P_0\{T_{nk}^o \geq a_{L-j+1}\},\ j = 1,\ldots,L.$$

Then $0 \leq \tau_1 < \tau_2 < \cdots < \tau_L \leq 1$. Also, let us write

(5.10) $$Y_k(\tau_j) = I(T_{nk}^o \geq a_{L-j+1}),\ j = 1,\ldots,L,\ k = 1,\ldots,K.$$

Further, let

(5.11) $$W_K(\tau_j) = \sum_{k=1}^{K} Y_k(\tau_j),\ j = 1,\ldots,L.$$

Also, let $J = \max\{j : 1 \leq j \leq L;\ \tau_j \leq \eta\}$ for some pre-assigned $\eta > 0$. Basically, we would like to pursue the distributional features of the partial sequence $\{W_K(\tau_j),\ j \leq J\}$, and incorporate Theorem 2. Note that in this way, we avoid the conventional assumption of a continuous null distribution of the coordinate-wise test statistics. Of course, if $n$ is adequately large, the assumption of a uniform distribution of the $p$-values (under the null hypothesis) would be reasonable. For example, if we have a three sample situation with $n_1 = n_2 = n_3 = 4$ then $L = (12)!/(4!)^3 = 34,650$ so that we could choose $J = 1$ and use the Poisson approximation. It is also possible to choose $J = 2$ with an appropriate cut-off point and still stick to a FWER around 0.05. In any case, under alternatives (of stochastic ordering) the distribution of the $T_{nk}^o$ will be tilted towards the right, still confined to the interval $(-1, 1)$, and hence, their centering would be shifted to the right of the origin with a negatively skewed distribution.

Corresponding to the known points $\tau_1 < \cdots < \tau_J$, let us consider the partial process $W_K(\tau_j), j = 1,\ldots,J$, as defined above. Also, let us choose a set of nonnegative integers $r_1 \leq \cdots \leq r_J$ in such a way that

(5.12) $$P_0\{W_K(\tau_j) > r_j, \text{for some } j \leq J\} = \alpha,$$

where $\alpha$ may not be exactly equal to a specified level (such as 0.05) but can be approximated very well through the above Poisson process result. If we let

(5.13) $$A_j = [\ W_K(\tau_j) > r_j\ ],\ j = 1,\ldots,J,$$

then (5.12) can be written as $P\{\ \bigcup_{j \leq J} A_j\ \}$, so that by the Bonferroni inequality,

(5.14) $$\begin{aligned}P\{\bigcup_{j \leq J} A_j\} &= \sum_{j \leq J} P\{A_j\} - \sum_{1 \leq j < j' \leq J} P\{A_j A_{j'}\} \\ &\quad + \sum_{1 \leq j < k < l \leq J} P\{A_j A_k A_l\} + \cdots + (-1)^K P\{A_1 \cdots A_K\}.\end{aligned}$$

Next, we use the Poisson approximation to each $P\{A_j\}$ wherein we use the following:

(5.15) $$P\{A_j\} \sim e^{-\nu_{nj}}\{\sum_{r > r_j} \nu_{nj}^r/r!\},\ j \geq 1.$$

Further, note that $W_K(\tau_j)$ is a nondecreasing (step) function in $j$ so that using the Markov property and Theorem 2 we may evaluate $P\{A_j A_{j'}\}$. Actually, we write $P\{A_j A_k\} = P\{A_j\} \cdot P\{A_k|A_j\}$, for $k > j$, and use Theorem 2 to approximate the conditional probability by $P\{Z_k > r_k | Z_j > r_j\}$ where $r_k \geq r_j, \forall j < k$. Also, typically terms involving more than 2 events $(A_j)$ will be small and can usually be neglected. Nevertheless, even if they are not small, the Markov property embedded in Theorem 2 can be used to provide a good approximation. Alternatively, we may write $P\{\bigcup_{j \leq J} A_j\} = 1 - P\{\bigcap_{j \leq J} A_j^c\}$ and using Theorem 2, write $P\{\bigcap_{j \leq J} A_j^c\}$ as a $J$-tuple sum over Poisson distributional probabilities. For small $r_j, j \leq J$, as is typically the case, this computation does not appear to be a formidable task.

Led by these findings, let us now consider the following testing procedure:

Compute the $W_K(\tau_j)$, $j \leq J$ as above. If $W_K(\tau_j) \leq m_j$, $\forall j \leq J$, accept the null hypothesis that there is no DG. On the other hand, if $W_K(\tau_j)$ is greater than $m_j$ for at least one $j \leq J$, then reject the null hypothesis that all the genes are NDG, and proceed to detect those genes $k \in \mathcal{K}$ as DG where

(5.16) $$\mathcal{K} = \{k \in \{1, \ldots, K\} : Y_k(\tau_j) = 1, \text{ for some } j \leq J\}.$$

Note that if for some $k$, $Y_k(\tau_j) = 1$ for some $j \leq J$, then $Y_k(\tau_{j'}) = 1$, $\forall j' \geq j$. Further, note that $\mathcal{K}$ is a stochastic subset of $\{1, \ldots, K\}$, and $R =$ cardinality of $\mathcal{K}$ is a (nonnegative) integer valued random variable. The overall significance level of this testing procedure is well approximated by the preassigned level $\alpha$.

Let us denote the following exclusive events by

(5.17) $$B_1 = A_1; \quad B_j = A_1^c \cdots A_{j-1}^c A_j, \ j \leq J.$$

Then, by definition, $A_j = \bigcap_{j \leq J} B_j$. With the same notation as in (4.7)—(4.11), we study the other measures (viz., PCER, PFER, FDR and $^p$FDR). Towards this, we consider the nonnull situation where $K_0$ are NDG and $K_1 = K - K_0$ are DG. To handle the distribution of $R$, the total number of rejections, we let

(5.18) $$\tau_j^* = (K_0 \tau_j + K_1 \beta_j)/K = \tau_j + (K_1/K)(\beta_j - \tau_j), \ j \geq 1,$$

where

(5.19) $$\beta_j = K_1^{-1} \sum_{k \in DG} P\{T_{nk}^o \geq a_{L-j+1} | k \in \{1, \ldots, K\} - \mathcal{K}_0\},$$

for $j = 1, \ldots, J$. Note by arguments similar to those in Sections 3 and 4, $\beta_j \gg \tau_j, \forall j \leq J$. We may write

(5.20) $$E(m_1) = \sum_{j \leq J} E(m_1 I(B_j)).$$

Next note that the events $B_j, j \leq J$, depend on the partial process $W_K(\tau_j), j \leq J$ and are thereby governed by Theorem 2 with $\nu_{nj} = K t_j^*, j \leq J$. On the other hand, the distribution of $m_1$ is governed by the process $W_K^o(\tau_j), j \leq J$, where the drift function for $W_K^o(\tau_j))$ is $\nu_{nj}^o = K_0 \tau_j, j \leq J$. Using Theorem 2 and the reproductive property of the Poisson distribution, we may well approximate the (conditional) distribution of $m_1$, given $R$, by a binomial law with parameters $(R, K_0 \tau_j/(K_0 \tau_j + K_1 \beta_j))$ whenever $B_j$ holds. Thus, we are able to provide a good approximation to the PFER by writing

(5.21) $$E(m_1) = \sum_{j \leq J} E\{E(m_1/R \mid R, B_j) R I(B_j)\}.$$

If $J = 1$, the conditional binomial law directly applies and we have the approximation

$$(5.22) \qquad \frac{K_0 t_1}{K_0 t_1 + K_1 t_1^*} . E(RI(R > r_1)) = K_0 t_1 P\{R \geq r_1\},$$

where the last step follows from the fact that for a Poisson variable $X$ with parameter $np$, $E(XI(X > r)) = npP\{X \geq r\}$. For $J \geq 2$, we have to apply the conditional binomial law under the sets $B_j$, followed by the distribution of $R$ over the sets $B_j$, and this can be done by repeated quadrature procedures. Numerical studies have thereby good scope.

By construction, rejection of the null hypothesis $H_0$ entails that $R > r_1$ and may even be greater than $r_1$ if $B_j$ pertains for some $j \geq 1$. As such, we do not have any problem in applying the original definition of FDR (in (4.10)). We write

$$(5.23) \qquad \text{FDR} = E(Q) = \sum_{j \leq J} E(QI(B_j)) = \sum_{j \leq J} E\{E(Q|R \in B_j)I(B_j)\}$$

and use the conditional binomial law for each term in the right hand side. Detailed numerical study is planned for a future communication.

We conclude this section with some pertinent remarks and observations. First, the use of the Chen–Stein theorem in a multi-state context can be done under fairly mild regularity conditions regarding the dependence of the genes. Secondly, by our choice of the $r_j, j \leq J$ and allowing possibly $J \geq 1$, we are not only in a position to allow more flexibility in the choice of statistical inference procedures but also to enforce the rejection of null hypothesis under a more structured setup. This allows us to study the FDR, etc., under more diverse setups. Further, using Kendall's tau statistic for each gene separately, we are in a position to allow heterogeneity of the gene expressions across the $K$ genes in a completely arbitrary manner, while under the null hypothesis, the distribution of the $T_{nk}^o, k = 1, \ldots, K$ being completely known provides an easy access to the incorporation of the Chen–Stein theorem. Finally, instead of using Kendall's tau statistic (coordinate-wise), it might be attractive to use more general rank statistics [17]. Though the distribution-free aspect holds under the null hypothesis, such distributions are more complex to evaluate and the associated Poisson processes have more complex drift functions. Further, such linear rank statistics involve some design variables which assume more structure on the $F_{ik}, k = 1, \ldots, K$, not necessary with the use of Kendall's tau.

**Acknowledgments.** The author is grateful to the reviewers for their critical reading of the manuscript and most helpful comments. Thanks are also due to Dr. Moonsu Kang and Sunil Suchandran for providing the Figure in the text.